\documentclass[a4paper,12pt]{article}
\usepackage{amsmath, amsfonts, amssymb}
\usepackage[english]{babel}
\makeindex
\begin{document}
\begin{center}

\textbf{\large Hyperspaces of superparacompact spaces and continuous maps}

\medskip
\textbf{A.A.Zaitov$^{1}$, D.I.Jumaev$^{2}$}

{$^{1}$ Tashkent institute of architecture and civil engineering,}

{adilbek\_zaitov@mail.ru},

{$^{2}$ Tashkent institute of architecture and civil engineering,}

{d-a-v-ron@mail.ru}

\end{center}

\sloppy
\begin{abstract}
In the present paper we establish that the space $\exp_\beta X$ of compact subsets of a Tychonoff space $X$ is superparacompact iff $X$ is so. Further, we prove the Tychonoff map $\exp_{\beta} f:\ \exp_{\beta} X\rightarrow \exp_{\beta} Y$ is superparacompact iff a given map $f:\ X\rightarrow Y$ is superparacompact.\\

{\bf Key words: } superparacompact space, finite-component cover, hyperspace, Tychonoff map.
\end{abstract}
\smallskip

\textbf{\large 0. Introduction}

In the present article under space we mean a topological $T_1$-space, under compact a Hausdorff compact space and under map a continuous map.

A collection $\omega$ of subsets of a set $X$ is said to be {\it star-countable} (respectively, {\it star-finite}) if each element of $\omega$ intersects at most countable (respectively, finite) set of elements of $\omega$.
A finite sequence of subsets  $M_{0}, ...,M_{s}$ of a set $X$ is a {\it chain} connecting sets $M_{0}$ and $M_{s}$, if $M_{i-1}\cap M_{i}\neq \varnothing$ for $i=1,...,s$. A collection $\omega$ of subsets of a set $X$ is said to be {\it connected} if for any pair of sets $M$, $M^{'}\subset X$ there exists a chain in $\omega$ connecting sets $M$ and $M^{'}$. The maximal connected subcollections of $\omega$ are called {\it components} of $\omega$.

For a system $\omega=\{O_{\alpha}: a\in A\}$ of subsets of a space $X$ we suppose $[\omega]=[\omega]_{X}=\{[O_{\alpha}]_{X}: a\in A\}$. For a space $X$, its some subspace $W$ and a set $B\subset X\setminus W$ (respectively, a point $x\in X\setminus W$) they say [2] that an open cover $\lambda$ of the space $W$ pricks out the set $B$ (respectively, the point $x$) in $X$ if $B\cap(\cup[\lambda]_{X})=\varnothing$ (respectively, $x\notin \cup [\lambda]_{X}$).

A star-finite open cover of a spaces $X$ is said to be {\it a finite-component cover} if the number of elements of each component is finite.

A collection $\omega$ of subsets of a set $X$ {\it refines} a collection $\Omega$ of subsets of $X$ if for each element $A\in \omega$ there is an element $B\in \Omega$ such that $A\subset B$. They also say that $\omega$ is {\it a refinement} of $\Omega$.

A space $X$ is said to be \textit{superparacompact} if every open cover of $X$ has  a finite-component cover which refines it [2].

Note that any compact is superparacompact, and any superparacompact space is strongly paracompact. Infinite discrete space is superparacompact, but it is not compact. Real line is strongly paracompact, but it is not superparacompact.

The following criterion plays a key role in investigation the class of superparacompact spaces.

\textbf{Theorem 1}[3]. A Tychonoff space $X$ is superparacompact iff for every closed set $F$ in $\beta X$ lying in the growth $\beta X\setminus X$ there exists a finite-component cover $\lambda$ of $X$ pricking out $F$ in $\beta X$ (i. e. $F\cap(\cup[\lambda]_{\beta X})=\varnothing$).

D.Buhagiar and T.Miwa offered the following criterion of superparacompactness.

\textbf{Theorem 2}[4]. A Tychonoff space $X$ is superparacompact iff for every closed set $F$ in perfect compactification $bX$ lying in the growth $bX\setminus X$ there is a finite-component cover $\lambda$ of $X$ pricking out $F$ in $bX$  (i. e. $F\cap(\cup[\lambda]_{bX})=\varnothing$).

For a topological space $X$ and its subset $A$ a set
$Fr_{X}A=[A]_{X}\cap[X\setminus A]_{X}=[A]_{X}\setminus Int_{X}A$
is called a boundary of $A$.

Let $vX$ be a compact extension of a Tychonoff space $X$. If $H\subset X$ is an open set in $X$, then by $O(H)$ (or by $O_{vX}(H)$) we denote a maximal open set in $vX$ satisfying $O_{vX}(H)\cap X=H$. It is easy to see that
$$O_{vX}(H)=\bigcup_{\substack{\Gamma\in\tau_{vX}, \\ \Gamma\cap X=H}}\Gamma,$$
where $\tau_{v X}$ is topology of the space $vX$.

A compactification $vX$ of a Tychonoff space $X$ is called {\it perfect with respect to an open set} $H$ in $X$, if the equality $[Fr_{X}H]_{vX}=Fr_{vX}O_{vX}(X)$ holds. If $vX$ is  perfect with respect to every open set in $X$, then it is called {\it a perfect compactification} of the space $X$ ([1], P. 232). A compactification $vX$ of a space $X$ is perfect iff for any two disjoint open sets $U_{1}$ and $U_{2}$ in $X$ the equality $O(U_{1}\bigcup U_{2})=O(U_{1})\bigcup O(U_{2})$ is executed. The Stone-C\v{e}ch compactification  $\beta X$ of $X$ is perfect. The equality $O(U_{1}\bigcup U_{2})=O(U_{1})\bigcup O(U_{2})$ is satisfied for every pair of open sets $U_{1}$ and $U_{2}$ in $X$ iff $X$ is normal, and the compactification $vX$ coincides with the Stone-C\v{e}ch compactification $\beta X$, i. e. $vX\cong\beta X$.

Let $X$ be a space. By $\exp X$ we denote a set of all nonempty closed subsets of $X$. A family of sets of the view
$$O\langle U_{1},...,U_{n}\rangle=\{F\in \exp X:\ F\subset \bigcup_{i=1}^n U_{n}, F\cap U_{1}\neq\varnothing,..., F\cap U_{n}\neq \varnothing\}$$
forms a base of a topology on $\exp X$, where $U_{1},\dots,U_{n}$ are open nonempty sets in $X$. This topology is called {\it the Vietoris topology}. A space $\exp X$ equipped with Vietoris topology is called {\it hyperspace} of $X$. For a compact $X$ its hyperspace $\exp X$ is also a compact.

Note for any space $X$ it is well known that $$\left[O\langle U_{1},...,U_{n}\rangle\right]_{\exp X}=O\left\langle [U_{1}]_{X},...,[U_{n}]_{X}\right\rangle.$$

Let $f:X\rightarrow Y$ be continuous map of compacts, $F\in \exp X$. We put
$$(\exp f)(F)=f(F).$$
This equality defines a map $\exp f: \exp X\rightarrow \exp Y$. For a continuous $f$ the map $\exp f$ is continuous. Really, it follows from the formula
$$
(\exp f)^{-1}O\langle U_{1},..., U_{m}\rangle=O\langle f^{-1}(U_{1}),..., f^{-1}(U_{m})\rangle
$$
what one can check directly. Note that if $f:X\rightarrow Y$ is an epimorphism, then $\exp f$ is also an epimorphism.

For a Tychonoff space $X$ we put
$$\exp_{\beta} X=\{F\in\exp \beta X: F\subset X\}.$$

It is clear, that $\exp_{\beta}X\subset \exp X$. Consider the set $\exp_{\beta}X$ as a subspace of the space $\exp X$. For a Tychonoff spaces $X$ the space $\exp_{\beta} X$ is also a Tychonoff space with respect to the induced topology.

For a continuous map $f:X\rightarrow Y$ of Tychonoff spaces we put $$\exp_{\beta}f =\left. (\exp\beta f)\right|_{\exp_\beta X},$$ where $\beta f: \beta X\rightarrow \beta Y$ is the  Stone-C\v{e}ch compactification of $f$ (it is unique).\\

\textbf{\large 1. Hyperspace of superparacompact spaces}

It is well known that for a Tychonoff space $X$ the set $\exp_{\beta} X$ is everywhere dense in $\exp \beta X$, i. e. $\exp \beta X$ is a compactification of the space $\exp_{\beta}X$. We claim $\exp \beta X$ is a perfect compactification. At first we will prove the following technical statement.

\textbf{Lemma 1.} Let $\gamma X$ be a compact extension of a space $X$ and, $V$ and $W$ be disjoint open sets in $\gamma X$. Let $V^{X}=X\cap V$ and $W^{X}=X\cap W$. Then the following equality is true:
$$[X\setminus V^{X}]_{\gamma X}\cap [X\setminus W^{X}]_{\gamma X}=[X\setminus(V^{X}\cup W^{X})]_{\gamma X}.$$

\textbf{Proof.} It is clear that $[X\setminus V^{X}]_{\gamma X}\cap [X\setminus W^{X}]_{\gamma X}\supset [X\setminus(V^{X}\cup W^{X})]_{\gamma X}$. Let $x\in [X\setminus V^{X}]_{\gamma X}\cap [X\setminus W^{X}]_{\gamma X}$. Then each open neighborhood $Ox$ in $\gamma X$ of $x$ intersects with the sets $X\setminus V^{X}$ and $X\setminus W^{X}$. Hence, $Ox\not\subset V^X$ and $Ox\not\subset W^X$. Therefore, since $V^X\cap W^X=\varnothing$, we have $Ox\not\subset V^X\cup W^X$, i. e. $Ox\cap X\setminus (V^{X}\cup W^{X})\neq \varnothing$. By virtue of arbitrariness of the neighbourhood $Ox$ we conclude that $x\in [X\setminus(V^{X}\cup W^{X})]_{\gamma X}.$ Lemma 1 is proved.

\textbf{Theorem 3.} For a Tychonoff space $X$ the space $\exp \beta X$ is a perfect compactification of the space $\exp_{\beta} X$.

\textbf{Proof.} It is enough to consider basic open sets. Let $U_{1}$ and $U_{2}$ be disjoint open sets in $X$. Since $\beta X$ is perfect compactification of $X$ we have $O_{\beta X}(U_{1}\cup U_{2})=O_{\beta X}(U_{1})\cup O_{\beta X}(U_{2})$. Consider open sets
$$O\langle U_{i}\rangle=\{F: F\in \exp_{\beta}X, F\subset U_{i}\},\qquad i=1,\ 2 $$
in $\exp_{\beta} X$. It is clear, that $O\langle U_{1}\rangle \cap O\langle U_{2}\rangle=\varnothing$. We will show that
$$O_{\exp \beta X}(O\langle U_{1}\rangle \cup O\langle U_{2}\rangle)=O_{\exp \beta X}(O\langle U_{1}\rangle)\cup O_{\exp \beta X}(O\langle U_{2}\rangle).$$
The inclusion $\supset$ follows from the definition of the set $O(H)$ (see [1], P. 234). That is why it is enough to show the inverse inclusion. Let $\Phi\subset \beta X$ be a closed set such that $\Phi\notin O_{\exp \beta X}(O\langle U_{1}\rangle)\cup O_{\exp \beta X}(O\langle U_{2}\rangle)$. Then $\Phi\in \exp \beta X \setminus O_{\exp \beta X}(O\langle U_{i}\rangle)$, $i=1,\ 2$. From [1] (see, P. 234) we have
$$\exp \beta X\setminus O_{\exp \beta X}(O\langle U_{i}\rangle)=[\exp_{\beta} X\setminus O\langle U_{i}\rangle]_{\exp \beta X}, \qquad  i=1,\ 2.$$
Hence $\Phi\in [\exp_{\beta} X\setminus O\langle U_{i}\rangle]_{\exp \beta X}$, $i=1,\ 2$. Since $O\langle U_{1}\rangle \cap O\langle U_{2}\rangle=\varnothing$ by Lemma 1 we have
$$[\exp_{\beta}X\setminus O\langle U_{1}\rangle]_{\exp \beta X}\cap [\exp_{\beta}X\setminus O\langle U_{2}\rangle]_{\exp \beta X}=[\exp_{\beta}X\setminus O(\langle U_{1}\rangle \cup O\langle U_{2}\rangle )]_{\exp \beta X}.$$
Therefore, $\Phi \in [\exp_{\beta}X\setminus O_{\exp \beta X}(O\langle U_{1}\rangle \cup O\langle U_{2}\rangle )]_{\exp \beta X}$, what is equivalent $\Phi \in \exp \beta X\setminus O_{\exp \beta X}(\langle U_{1}\rangle \cup \langle U_{2}\rangle )$  (see [1], P. 234). In other words, $\Phi \notin O_{\exp \beta X}(\langle U_{1}\rangle \cup \langle U_{2}\rangle )$. Thus, we have established that inclusion
$O_{\exp \beta X}(\langle U_{1}\rangle \cup \langle U_{2}\rangle )\subset O_{\exp \beta X}(O\langle U_{1}\rangle\cup O_{\exp \beta X}(O\langle U_{2}\rangle $ is also fair. Theorem 3 is proved.

\textbf{Lemma 2.} Let $U_1,\dots,\ U_n;\ V_1,\dots,\ V_m$ be open subsets of a space $X$. Then $O\langle U_1,\dots,\ U_n\rangle\cap O\langle V_1,\dots,\ V_m\rangle\neq\varnothing$ iff for each $i \in \{1,\dots, n\}$ and for each $j\in \{1,\dots, m\}$ there exists, respectively $j(i)\in \{1, \dots, m\}$ and $i(j) \in \{1,\dots, n\}$, such that $U_i\cap V_{j(i)}\neq \varnothing$ and $U_{i(j)}\cap V_j\neq \varnothing$.

\textbf{Proof.} Assume that for every $i \in \{1,\dots, n\}$ there exists $j(i)\in \{1,\dots, m\}$ such that $U_{i(j)}\cap V_{j(i)}\neq \varnothing$ and for every $j\in \{1,\dots, m\}$ there exists $i(j) \in \{1,\dots, n\}$ such that $U_{i(j)}\cap V_j\neq \varnothing$. For any pair $(i,\ j)\in \{1,\dots, n\}\times \{1,\dots, m\}$ for which $U_i\cap V_j\neq \varnothing$, choose a point $x_{ij}\in U_i\cap V_j$ and make a closed set $F$ consisting of these points. Then $F\subset \bigcup\limits_{i=1}^{n}U_i$ and $F\subset \bigcup\limits_{j=1}^{m}V_j$. Besides, $F\cap U_i\neq\varnothing$, $i=1,\dots,\ n$, and $F\cap V_j\neq\varnothing$, $j=1,\dots,\ m$. Therefore, $F\in O\langle U_1,\dots,\ U_n\rangle\cap O\langle V_1,\dots,\ V_m\rangle$.

Suppose there exists $i_0\in \{1,\dots, n\}$ such that $U_{i_0}\cap V_{j} = \varnothing$ for all $j\in \{1,\dots, m\}$. Then $U_{i_0}\cap \bigcup\limits_{j=1}^{m}V_j= \varnothing$ and for each $F\in O\langle U_1,\ ...,\ U_n\rangle$ we have $F\not\subset \bigcup\limits_{j=1}^{m}V_j$. Hence, $F\notin O\langle V_1,\ ...,\ V_m\rangle$. Similarly, every $\Gamma\in O\langle V_1,\ ...,\ V_m\rangle$ lies in $\bigcup\limits_{j=1}^{m}V_j$ what implies $\Gamma\cap U_{i_0}=\varnothing$. From here $\Gamma\notin O\langle U_1,\ ...,\ U_n\rangle$. Thus, $O\langle U_1,\ ...,\ U_n\rangle\cap O\langle V_1,\ ...,\ V_m\rangle=\varnothing$. Lemma 2 is proved.

\textbf{Lemma 3.}  Let $\upsilon$ be a finite-component cover of a Tychonoff space $X$. Then the family $\exp_\beta \upsilon =\{O\langle U_1,\dots,\ U_n\rangle: U_i\in \upsilon,\ i=1,\dots,\ n;\ n\in \mathbb{N}\}$ is a finite-component cover of the space $\exp_{\beta} X$.

\textbf{Proof.}  Let $O\langle G_1,\dots,\ G_k\rangle$ be an element of $\exp_{\beta} \upsilon$. Each $G_i\in \upsilon$ intersects with finite elements of $\upsilon$. Let $\left|\{\alpha:\ G_i\cap U_\alpha\neq\varnothing,\ U_\alpha\in \upsilon \}\right|=n_i$, $i=1,\ 2,\dots,\ k$. Denote $\gamma=\{G_i\cap U_j:\ G_i\cap U_j\neq \varnothing,\ i=1,\ 2,\dots,\ k,\ U_j\in \upsilon\}$. Then $|\gamma|\leq n_1+\dots+n_k$. Therefore, the set $O\langle G_1,\dots,\ G_k\rangle$ crosses no more then $n_1+\dots+n_k$ elements of $\exp_\beta \upsilon$. It means that the collection $\exp_\beta \upsilon$ is star-finite.

Let $F\in \exp_\beta X$. There is a subfamily $\upsilon_F \subset \upsilon$ such that $F\subset \bigcup\limits_{U\in\upsilon_F}U$. From a cover $\{F\cap U:\ U\in \upsilon_F,\ F\cap U\neq\varnothing\}$ of the compact $F$ it is possible to allocate a finite subcover $\{F\cap U_i:\ i=1,\dots, m\}$. We have $F\in O\langle U_1,\dots,\ U_m \rangle$. So, the family $\exp_\beta \upsilon$ is a cover of $\exp_{\beta}X$. On the other hand by the  definition of Vietoris topology the cover $\exp_\beta \upsilon$ is open. Thus, $\exp_\beta \upsilon$ is a star-finite open cover of $\exp_{\beta} X$.

We will show now that all components of the $\exp_{\beta} \upsilon$ are finite.

Let $M=O\langle G_1,\dots,\ G_s\rangle$ and $M^{'}=O\langle G_1^{'},\dots,\ G_t^{'}\rangle$ be arbitrary elements of $\exp_{\beta} \upsilon$. Further, let $\gamma_{G_{i}G_{j}^{'}}=\{U_{l}^{ij}:\ l=1,\ 2,\dots,\ n_{ij}\}$ be the maximal chain of $\upsilon$ connecting $G_{i}$ and $G_{j}^{'}$, $i=1,\ 2,\dots,\ s$, $j=1,\ 2,\dots,\ t$. By definition these sets satisfy the following properties:

$
(1)\qquad U_1^{ij}=G_i, \qquad i=1,\dots,s;\ j=1,\dots,t;
$

$
(2)\qquad U_{n_{ij}}^{ij}=G_j^{'}, \qquad i=1,\dots,s;\ j=1,\dots,t;
$

$
(3)\qquad U_l^{ij}\cap U_{l+1}^{ij}\neq \varnothing, \qquad l=1,\dots,n_{ij}-1;\ i=1,\dots,s;\ j=1,\dots, t.
$

If $s<t$ we have $O\langle G_1,\dots, G_s\rangle = O\langle U_1^{1j},\dots, U_1^{sj},\ U_1^{i_{1}(s+1)},\dots, U_1^{i_{t-s}t}\rangle$, where $j=1,\dots,t$ and $i_1,\dots,i_{t-s}\in \{1,\dots,s\}$. Further, $O\langle G_1^{'},\dots, G_t^{'}\rangle = O\langle U_{n_{i1}}^{i1},\dots, U_{n_{it}}^{it}\rangle$, $i=1,\dots, s$. Thus, the cover $\exp_{\beta} \upsilon$ has a chain connecting the given sets $M=O\langle G_1,\dots,\ G_s\rangle$ and $M^{'}=O\langle G_1^{'},\dots,\ G_t^{'}\rangle$. The case $s>t$ is analogously.

Now using Lemma 2 and calculating directly we find that each maximal chain of $\exp_{\beta} \upsilon$ connecting the sets $M=O\langle G_1,\dots,\ G_s\rangle$ and $M^{'}=O\langle G_1^{'},\dots,\ G_t^{'}\rangle$ has no more than $\prod\limits_{\substack{i=1,\\ j=1}}^{\substack{t\\ s}}n_{ij}$ elements. Thus, all components of $\exp_{\beta} \upsilon$ is finite. Lemma 3 is proved.

\textbf{Theorem 4.} For a Tychonoff space $X$ its hyperspace $\exp_{\beta}X$  is superparacompact iff $X$ is superparacompact.

\textbf{Proof.} As the superparacompactness is inherited to the closed subsets [2], the superparacompactness of $\exp_{\beta}X$  implies superparacompactness of the closed subset $X\subset \exp_{\beta}X$.

Let $\Omega$ be an open cover of $\exp_\beta X$. For each element $G\in \Omega$ there exists $O_G\langle U_1,\dots, U_n\rangle$ such that $O_G\langle U_1,\dots, U_n\rangle\subset G$, where $U_1,\dots, U_n$ are open sets in $X$. We can choose sets $G\in \Omega$ so that a collection of sets $O_G\langle U_1,\dots,\ U_n\rangle$ forms a cover of $\exp_\beta X$, what we denote by $\Omega^{'}$. It easy to see that a collection $\omega^{'}=\bigcup\limits_{O_G\langle U_1,\dots, U_n\rangle \in \Omega^{'}}\{U_1,\dots,\ U_n\}$ is an open cover of $X$. There exists a finite-component cover $\omega$ of $X$ which refines $\omega^{'}$. Then by Lemma 3 the collection
$$
\exp_\beta\omega=\{O\langle V_1,\dots,\ V_k\rangle:\ V_i\in \omega,\ i=1,\dots,\ n;\ n\in \mathbb{N}\}
$$
is a finite-component cover of $\exp_\beta X$ and it is refines $\Omega$. Theorem 4 is proved.\\

\textbf{\large 2. Superparacompactness of the map $\exp_{\beta} f$}

For a continuous map $f: (X, \tau_X)\rightarrow (Y, \tau_Y)$ and $O\in \tau_Y$ a preimage $f^{-1}O$ is called {\it a tube} (above $O$). Remind, a continuous map $f:\ X\rightarrow Y$ is called [2] $T_0$-{\it map}, if for each pair of distinct points $x$, $x^{'}\in X$, such that $f(x)=f(x^{'})$, at least one of these points has an open neighborhood in $X$ which does not contain another point. A continuous map $f:\ X\rightarrow Y$ is called {\it totally regular}, if for each point $x\in X$ and every closed set $F$ in $X$ not containing $x$ there exists an open neighborhood $O$ of $f(x)$ such that in the tube $f^{-1}O$ the sets $\{x\}$ and $F$ are functional separable. Totally regular $T_0$-map is said to be {\it Tychonoff map}.

Obviously, each continuous map $f:\ X\rightarrow Y$  of a Tychonoff space $X$ into a topological space $Y$ is a Tychonoff map. In this case owing to the set $\exp_{\beta} X$ is a Tychonoff space concerning to Vietoris topology for every Tychonoff space $X$, the map $\exp_{\beta} f:\ \exp_{\beta} X\rightarrow \exp_{\beta} Y$ is a Tychonoff map.

A continuous, closed map $f:X\rightarrow Y$ is said to be {\it compact} if the preimage $f^{-1}y$ of each point $y\in Y$ is compact. A continuous map $f: X\rightarrow Y$ is compact iff for each point $y\in Y$ and every cover $\omega$ of the fibre $f^{-1}y$, consisting of open sets in $X$, there is an open neighbourhood $O$ of $y$ in $Y$ such that the tube $f^{-1}O$ can be covered with a finite subfamily of $\omega$.

A compact map $bf:\ b_f X\rightarrow Y $ is said to be {\it a compactification} of a continuous map $f:X\rightarrow Y$ if $X$ is everywhere dense in $b_f X$ and $bf|_{X}=f$. On the set of all compactifications of the map $f$ it is possible to introduce a partial order: for the compactifications $b_{1}f:\ b_{1f}X\rightarrow Y$ and $b_{2}f:\ b_{2f}X\rightarrow Y$ of $f$ we put $b_{1}f\leq b_{2}f$ if there is a natural map of $b_{2f}X$ onto $b_{1f}X$. B. A. Pasynkov showed that for each Tychonoff map $f:X\rightarrow Y$ there exists its maximal compactification $g:Z\rightarrow Y$, which he denoted by $\beta f$, and the space $Z$ where this maximal compactification defines by $\beta_{f}X$. To within homeomorphism for a given Tychonoff map $f$ its maximal compactification $\beta f$ is unique.

A Tychonoff map $f:X\rightarrow Y$ is said to be {\it superparacompact}, if for every closed set $F$ in $\beta _{f}X$ lying in the growth $\beta_{f}X\setminus X$  there exists a finite-component cover $\lambda$ of $X$ pricking out $F$ in $\beta _{f}X$ (i. e. $F\cap(\cup[\lambda]_{\beta_f X})=\varnothing$) [3].

It is easy to see that one can determine superparacompactness of a map as so: a map $f:X\rightarrow Y$ is superparacompact if for each $y\in Y$ and every open cover $\Upsilon$ of $f^{-1}y$ in $X$ there exists an open neighbourhood $O$ of $y$ in $Y$ such that $\Upsilon$ has a finite-component cover $\upsilon$ of $f^{-1}O$ in $X$ which refines $\Upsilon$.

\textbf{Definition 1.} A compactification $bf: b_f X \rightarrow Y$ of a Tychonoff map $f: X\rightarrow Y$ is said to be {\it perfect compactification} of $f$ if for each point $y\in Y$ and for every disjoint open sets $U_{1}$ and $U_{2}$ in $X$ there exists an open neighbourhood $O\subset Y$ of $y$ such that the equality $$O_{b_f X}(U_{1}\cup U_{2})\cap bf^{-1}O=\left(O_{b_f X}(U_{1})\cup O_{b_f X}(U_{2})\right)\cap bf^{-1}O$$
holds.

Let $f:\ X\rightarrow Y$ be a continuous map of a Tychonoff space $X$ into a space $Y$. It is well known there exists a compactification $vX$ of $X$ such that $f$ has a continuous extension $vf:\ vX\rightarrow Y$ on $vX$. It is clear, $vf$ is a perfect compactification of $f$.

The following result is an analog of Theorem 2 for a case of maps.

\textbf{Theorem 5}. Let $bf:\ b_f X\rightarrow Y$ be a perfect compactification of a Tychonoff map $f:X\ \rightarrow Y$. The map $f$ is superparacompact iff for every closed set $F$ in $b_f X$ lying in the growth $b_f X\setminus X$ there exists a finite-component cover $\lambda$ of $X$ pricking out set $F$ in $b_f X$.

\textbf{Proof} is carried out similar to the proof of Theorem 1.1 $\Pi$ from [2].

Evidently a restriction $f|_{\Phi}:\ \Phi\rightarrow Y$ of a superparacompact map $f:\ X\rightarrow Y$ on the closed subset $\Phi\subset X$ is a superparacompact map.

The following result is a variant of Theorem 3 for a case of maps.

\textbf{Theorem 6.} Let $f:\ X\rightarrow Y$ be a Tychonoff map. Then the map $\exp_{\beta} \beta f : \exp_{\beta} \beta_f X\rightarrow \exp_{\beta} Y$ is a perfect compactification of $\exp_{\beta} f:\exp_{\beta} X\ \rightarrow \exp_{\beta} Y$.

\textbf{Proof} is similar to the proof of Theorem 3. Here the equality
$$
(\exp_{\beta} \beta f)^{-1}O\langle U_{1},..., U_{m}\rangle=O\langle \beta f^{-1}(U_{1}),..., \beta f^{-1}(U_{m})\rangle
$$
is used.

The following statement is the main result of this section.

\textbf{Theorem 7.} The Tychonoff map $\exp_{\beta} f:\ \exp_{\beta} X\rightarrow \exp_{\beta} Y$ is superparacompact iff a map $f:\ X\rightarrow Y$ is superparacompact.

\textbf{Proof.} Let $\exp_{\beta} f:\ \exp_{\beta} X\rightarrow \exp_{\beta} Y$ be a superparacompact map. It implies that $f:\ X\rightarrow Y$ is a superparacompact map since $X\cong\exp_1 X$ is closed set in $\exp_\beta X$.

Let now $f:\ X\rightarrow Y$ be a superparacompact map. Consider arbitrary $\Gamma\in\exp_\beta Y$ and an open cover $\Omega$ of $(\exp_\beta f)^{-1}(\Gamma)=\{F\in \exp_\beta X:\ f(F)=\Gamma\}$ in $\exp_{\beta} X$. For each element $G\in \Omega$ there exists $O_G\langle U_1,\dots, U_n\rangle$ such that $O_G\langle U_1,\dots, U_n\rangle\subset G$, where $U_1,\dots, U_n$ are open sets in $X$. We can choose sets $G\in \Omega$ so that a collection of sets $O_G\langle U_1,\dots,\ U_n\rangle$ forms a cover of $(\exp_\beta f)^{-1}(\Gamma)$, what we denote by $\Omega^{'}$. It easy to see that a collection $\omega^{'}=\bigcup\limits_{O_G\langle U_1,\dots, U_n\rangle \in \Omega^{'}}\{U_1,\dots,\ U_n\}$ is an open cover of $f^{-1}\Gamma$ in $X$. For each $y\in \Gamma$ the collection $\omega_{y}=\{U\cap f^{-1}y:\ U\in\omega^{'}\}$ is an open cover of $f^{-1}y$ in $X$. There exists an open neighbourhood $O_y$ of $y$ in $Y$ such that $\omega_{y}$ has a finite-component cover $\omega_{y}^{'}$ of $f^{-1}O_y$ in $X$ which refines $\omega_{y}$. Gather such $O_y$ and construct an open cover $\{O_y: y\in \Gamma\}$ of $\Gamma$ in $Y$. Since $\Gamma$ is a compact there exists a finite open subcover $\gamma=\{O_{y_1}, \dots, O_{y_n}\}$ in $Y$, which covers $\Gamma$ . Put $\omega=\bigcup\limits_{O_{y_i}\in \gamma}\omega_{y_i}^{'}$. Then $\omega$ is an open cover of $f^{-1}\bigcup\limits_{U\in \omega}U$ in $X$. By construction $\omega$ is a finite-component cover, and it refines $\omega^{'}$. Hence, $\exp_\beta \omega$ is a finite-component cover of $(\exp_\beta f)^{-1}O\langle O_{y_1}, \dots, O_{y_n}\rangle=\left\langle f^{-1}O_{y_1}, \dots, f^{-1}O_{y_n}\right\rangle$ in $\exp_{\beta} X$ and it refines $\Omega$.

So, for each $\Gamma \in \exp_\beta Y$ and every open cover $\Omega$ of $(\exp_\beta f)^{-1}\Gamma$ in $\exp_\beta X$ there exists an open neighbourhood $O\langle O_{y_1},\dots, O_{y_n}\rangle$ of $\Gamma$ in $\exp_\beta Y$ such that $\Omega$ has a finite-component cover $\exp_\beta \omega$ of $(\exp_\beta f)^{-1}O\langle O_{y_1}, \dots, O_{y_n}\rangle$ in $\exp_\beta X$ which refines $\Omega$. Thus, the  map $\exp_{\beta} f:\ \exp_{\beta} X\rightarrow \exp_{\beta} Y$ is superparacompact. Theorem 7 is proved.

\textbf{Corollary 1.} Let $f:\ X\rightarrow Y$ be a superparacompact map and $\Phi$ be a closed set in $\exp_{\beta} \beta_f X$ such that $\Phi\subset \exp_{\beta} \beta_f X \setminus \exp_{\beta} X$. Then there exists a finite-component cover $\Omega$ of $\exp_\beta X$ pricking out $\Phi$ in $\exp_\beta\beta_{f}X$ (i. e. $\Phi\cap(\cup[\Omega]_{\exp_\beta{\beta_f X}})=\varnothing$).

\textbf{Corollary 2.} The functor $\exp_\beta$ lifts onto category of superparacompact spaces and their continuous maps.

\textbf{Acknowledgments.} \textit{The authors would like to acknowledge the comprehensive backing of professor Sh.Ayupov.}

\begin{center}\textbf{References}\end{center}

1. A. V. Arkhangelsky, V. I. Ponomarev. Fundamentals of the general topology: Problems and Exercises. -- D.Reidel Publishing Company. 1983. -- 415 pp (Originally published as: {\it Osnovy Obsheii Topologii v Zadachakh i Uprajneniyakh}, by A. V. Arkhangelsky, V. I. Ponomarev, Izdatel'stvo `Nauka' Moscow, 1974).

2. D. K. Musayev, B. A. Pasynkov. On  compactness  and  completeness  properties  of  topological  spaces and continuous maps. -- Tashkent: `Fan'. 1994. -- 124 pp.

3. D. K. Musayev. On  compactness  and  completeness  properties  of  topological  spaces and continuous maps. -- Tashkent: `NisoPoligraf'. 2011. -- 216 pp.

4. D. Buhagiar, T. Miwa. On Superparacompact and Lindelof GO-Spaces. //Houston Journal of Mathematics, Vol. 24, No. 3, 1998. P. 443 -- 457.

5. R. Engelking. General Topology. -- Polish Scientific Publishers. Warszawa. -- 1977.

6. A. V. Zarelua. On a theorem of Hurewicz. //Mat. Sb. (N.S.).
1963. Vol 60(102). No. 1. P. 17–28.

7. E. G. Sklyarenko. Some problems of the theory of bicompact extensions. //Izv. Akad. Nauk SSSR. Ser. Matem. 1962. Vol. 26. No. 3. P. 427–452.

8. V. V. Fedorchuk, V. V. Filippov, "General Topology. Basic Constructions" (in Russian). -- Moscow. Fizmatlit. 2006.

9. Ren\'{e} Bartsch. Hyperspaces in topological Categories. //arXiv:1410.3137v2 [math.GN] September 4, 2018. P. 13.

10. D. K. Musaev, D. I. Jumaev. On Generalization of the Freudenthal's Theorem for Compact Irreducible Standard Polyhedric Representation for Superparacompact Complete Metrizable Spaces. //arXiv:1304.0541 [math.GN] April 2, 2013. P. 11.

11. Valentin Gutev. Hausdorff continuous sections. //J. Math. Soc. Japan. Vol. 66, No. 2 (2014) pp. 523–534.
doi: 10.2969/jmsj/06620523

\end{document}